\documentclass[preprint,review,11pt]{elsarticle}
\usepackage{amssymb,amsmath}
\usepackage[bookmarksnumbered,colorlinks, plainpages]{hyperref}
\newtheorem{theorem}{Theorem}[section]
\newtheorem{lemma}{Lemma}[section]
\newtheorem{corollary}{Corollary}[section]
\newtheorem{proposition}{Proposition}[section]
\usepackage{enumerate}

\newtheorem{remark}{Remark}[section]

\journal{arXiv}

\begin{document}
\begin{frontmatter}

\title{\textbf{A note on the proofs of generalized Radon inequality}}

\author {Yongtao Li\footnote{\textit{E-mail address: liyongtaosx@outlook.com}},
Xian-Ming Gu$^{\dag,}$\footnote{\textit{${\dag}$ Corresponding author. 
E-mail address: guxianming@live.cn}},
Jianci Xiao\footnote{\textit{E-mail address: jcxshaw@outlook.com}}}

\address{
\begin{center}
School of Mathematics and Statistics, Central South University,\\
Changsha, Hunan, 410083, P.R. China
\end{center}

\begin{center}
School of Mathematical Sciences, University of Electronic Science 
and Technology of China, Chengdu, Sichuan, 611731, P.R. China
\end{center}

\begin{center}
School of Mathematical Sciences, Zhejiang University, \\
Hangzhou, Zhejiang, 310027, P.R. China
\end{center}
}

\begin{abstract}
In this paper, we introduce and prove the generalizations of Radon 
inequality. The proofs in the paper unify and are simpler than those 
in former work. Meanwhile, we also find mathematical equivalences 
among the Bernoulli inequality, the weighted AM-GM inequality, the 
H\"{o}lder inequality, the weighted power mean inequality and the 
Minkowski inequality. Finally, a series of the applications are shown 
in this note.

\end{abstract}

\begin{keyword}
Bergstr\"{o}m inequality, Radon inequality, Weighted power mean inequality, 
Equivalence, H\"{older} inequality.

\end{keyword}

\end{frontmatter}
\section{Introduction}
\setcounter{equation}{0}
\renewcommand{\theequation}{1.\arabic{equation}}

The well-known  Bergstr\"{o}m inequality (see e.g. \cite{KYF,HB,RB}) says 
that if $x_k,y_k$ are real numbers and $y_k >0$ for all $1\le k\le n$, then
\begin{equation}
\frac{x_1^2}{y_1} + \frac{x^{2}_2}{y_2} + \cdots + \frac{x^{2}_n}{y_n} \ge 
\frac{\left(x_1 + x_2 + \cdots + x_n \right)^2}{y_1 + y_2 + \cdots + y_n}
\label{e1}
\end{equation}
and the equality holds if and only if $\frac{x_1}{y_1} = \frac{x_2}{y_2} 
= \cdots=\frac{x_n}{y_n}$.

Some generalizations of the inequality (\ref{e1}) can be found in \cite{SBJ,JFY}. 
Actually, the following Radon inequality \eqref{e2} is just a direct consequence: 
If $b_1,b_2,\ldots,b_n$ are positive real numbers and $a_1, a_2,\ldots,a_n$, 
$m$ are nonnegative real numbers, then
\begin{equation}
\frac{a^{m+1}_1}{b^{m}_1} + \frac{a^{m+1}_2}{b^{m}_2} + \cdots + \frac{a^{m+1}_n}
{b^{m}_n} \ge \frac{\left(a_1 + a_2 + \cdots + a_n \right)^{m+1}}{\left(b_1 + 
b_2 + \cdots + b_n \right)^m}.
\label{e2} 
\end{equation}
When $m=1$, \eqref{e2} reduces to \eqref{e1}. For more details about Radon inequality 
\eqref{e2}, the readers can refer to \cite[pp. 1351]{JR} and \cite{GJG, CBM,CM}. In 
fact, it is not hard to prove that (\ref{e1}) is equivalent to the Cauchy-Buniakovski-Schwarz 
inequality (see \cite[pp. 34-35, Theorem 1.6.1]{RJR}) stated as follows: if $a_1, 
\ldots, a_n, b_1, \ldots, b_n$ are nonnegative real numbers, then 
\begin{equation*} \sum\limits_{k=1}^n a_k \sum\limits_{k=1}^n b_k \ge  
\left( \sum\limits_{k=1}^n \sqrt{a_kb_k}\right)^2. 
\end{equation*}

In \cite[Theorem 1]{YKC}, Yang has given a generalization of Radon inequality as 
follows: if $a_1, a_2,\ldots,a_n$ are nonnegative real numbers and $b_1,b_2,
\ldots,b_n$ are positive real numbers,  then for $r\ge 0,s\ge 0$ and $r\ge 
s + 1$,
\begin{equation}
\frac{a^{r}_1}{b^{s}_1} + \frac{a^{r}_2}{b^{s}_2} + \cdots + \frac{a^{r}_n}
{b^{s}_n} \ge \frac{\left(a_1 + a_2 + \cdots + a_n\right)^r}{n^{r - s - 1}
\left(b_1 + b_2 + \cdots + b_n\right)^s}.
\label{ee5}
\end{equation}

The weighted power mean inequality (see \cite[pp. 111-112, Theorem 10.5]{ZC}, 
\cite[pp. 12-15]{GJG} and \cite{DJA} for details) is defined as follows: if 
$x_1, x_2,\ldots, x_n$ are nonnegative real numbers and $p_1, p_2, \ldots, 
p_n$ are positive real numbers, then for $r \ge s >0 $, we have
\begin{equation}
\left(\frac{p_1 x^{r}_1 + p_2 x^{r}_2 + \cdots + p_n x^{r}_n}{p_1 + p_2 + 
\cdots + p_n}\right)^{\frac{1}{r}} \ge \left(\frac{p_1 x^{s}_1 + p_2 x^{s
}_2 + \cdots + p_n x^{s}_n}{p_1 + p_2 + \cdots + p_n}\right)^{\frac{1}{s}}.
\label{e3}
\end{equation}

In this paper, we give three different cheaper proofs and some applications 
of generalized Radon inequality \eqref{ee5}, and then present equivalence 
relations between the weighted power mean inequality and Radon inequality. 
Furthermore, we summarize the equivalences among the weighted AM-GM inequality, 
the H\"{o}lder inequality, the weighted power mean inequality and the 
Minkovski inequality.

\section{Main results}
\setcounter{equation}{0}
\renewcommand{\theequation}{2.\arabic{equation}}

In this section, we first give three different methods for proving the generalized 
Radon inequality \eqref{ee5}. To read for convenience, the result obtained by Yang 
can be stated as the following theorem.

\begin{theorem} If $a_1, a_2,\ldots,a_n$ are nonnegative real numbers and $b_1,
b_2,\ldots,b_n$ are positive real numbers, then for $r\ge 0,s\ge 0$ and $r \ge 
s + 1$,
\begin{equation}
\frac{a_1^r}{b_1^s}+\frac{a_2^r}{b_2^s}+\cdots+\frac{a_n^r}{b_n^s} \ge \frac{\left(a_1+a_2+\cdots+a_n\right)^r}{n^{r-s-1}\left(b_1+b_2+\cdots+b_n\right)^s}.
\label{e5}
\end{equation}
\end{theorem}

\noindent
{\bf Proof 1.}
By Radon inequality (\ref{e2}), we have
\begin{eqnarray}\label{e6}\sum\limits_{k=1}^n \frac{a_k^r}{b_k^s} = \sum\limits_{k=1}^n \frac{\left(a_k^{\frac{r}{s+1}}\right)^{s+1}}{b_k^s} \ge \frac{\left( a_1^{\frac{r}{s+1}}+a_2^{\frac{r}{s+1}}+\cdots+a_n^{\frac{r}{s+1}} \right)^{s+1}}{ \left( b_1+b_2+\cdots+b_n \right)^s}.\end{eqnarray}
Note that $r\ge s+1\ge 1$, then $\frac{r}{s+1} -1\ge0$. Using Radon inequality again, we get that
\begin{eqnarray}\label{e7} \sum\limits_{k=1}^n a_k^{\frac{r}{s+1}} = \sum\limits_{k=1}^n \frac {a_k^{\frac{r}{s+1}}}{1^{\frac{r}{s+1}-1}}\ge\frac{\left(a_1+a_2+\cdots+a_n \right)^{\frac{r}{s+1}}}{\left(1+1+\cdots+1\right)^{\frac{r}{s+1}-1}} .\end{eqnarray}
According to inequalitis (\ref{e6}) and (\ref{e7}), we clearly have
\begin{eqnarray*} \frac{a_1^r}{b_1^s}+\frac{a_2^r}{b_2^s}+\cdots+\frac{a_n^r}{b_n^s} \ge \frac{\left(a_1+a_2+\cdots+a_n\right)^r}{n^{r-s-1}\left(b_1+b_2+\cdots+b_n\right)^s}.\end{eqnarray*}
Therefore, the desired result (\ref{e5}) is obtained.

\noindent
{\bf Proof 2.}
Let the concave  function $f$ : $(0, +\infty)$ $\rightarrow \mathbb{R}$ be $\ln x$. We observe that the weighted Jensen inequality: for $q_1, q_2, q_3 \in [0,1]$ with $q_1+q_2+q_3=1$ and positive real numbers $x_1,x_2,x_3$, then we have
\begin{eqnarray*} q_1f(x_1)+q_2f(x_2)+q_3f(x_3) \le f(q_1x_1+q_2x_2+q_3x_3),\end{eqnarray*}
and the equality holds if and only if $x_1=x_2=x_3$.
We denote $$U_n(a)=\left( \frac{a_1^r}{b_1^s}+\frac{a_2^r}{b_2^s}+\cdots+\frac{a_n^r}{b_n^s} \right)^{-1}$$ and $$H_n(b)=(b_1+b_2+\cdots+b_n)^{-1}.$$
Consider $x_1=\frac{a_k^r}{b_k^s}U_n(a), x_2=b_k H_n(b),x_3=\frac{1}{n}$ and $q_1=\frac{1}{r},q_2=\frac{s}{r},q_3=\frac{r-s-1}{r}$ (observe that $q_3 \ge 0$ from $r \ge s+1$). So we have
\begin{eqnarray*}\begin{aligned}&a_k (U_n(a))^{\frac{1}{r}} \cdot(H_n(b))^{\frac{s}{r}} \cdot \left(\frac{1}{n}\right)^{\frac{r-s-1}{r}} \\
&\le  \frac{1}{r} \cdot \frac{a_k^r}{b_k^s}U_n(a)+ \frac{s}{r}\cdot b_k H_n(b)+ \frac{r-s-1}{r}\cdot \frac{1}{n}.\end{aligned}\end{eqnarray*}
Summing up over $k$ $(k=1,2,\ldots, n)$, we obtain
\begin{eqnarray*}\begin{aligned}&\sum\limits_{k=1}^na_k (U_n(a))^{\frac{1}{r}} \cdot(H_n(b))^{\frac{s}{r}} \cdot \left(\frac{1}{n}\right)^{\frac{r-s-1}{r}}\\
&\le \sum\limits_{k=1}^n \left( \frac{1}{r}\cdot \frac{a_k^r}{b_k^s}U_n(a)+ \frac{s}{r}\cdot b_k H_n(b)+ \frac{r-s-1}{r} \cdot \frac{1}{n} \right)=1.\end{aligned}\end{eqnarray*}
The required inequality (\ref{e5}) follows.

For many numerical inequalities, the induction is some times a useful method used to establish a given statement for all natural numbers. We now give the third proof of Theorem 2.1 by mathematical induction. To state this proof clearly, let us start with the following lemma.

\begin{lemma}\label{lemma2.1} If $a_1, a_2, \ldots, a_n,b_1, b_2, \ldots, b_n$ are nonnegative real numbers and $\lambda_1$, $\lambda_2$,$\ldots$,$\lambda_n$ are nonnegative real numbers such that $\lambda_1+\lambda_2+\cdots+\lambda_n=1$, then
\begin{eqnarray}\label{e8}\prod\limits_{k=1}^n a_k^{\lambda_k}+\prod\limits_{k=1}^n b_k^{\lambda_k} \le \prod\limits_{k=1}^n \left(a_k+b_k \right)^{\lambda_k}.\end{eqnarray}\end{lemma}

\noindent
{\bf Proof of lemma 2.1. }According to the weighted AM-GM inequality, we have
\begin{eqnarray*}\label{e9}\prod\limits_{k=1}^n\left(\frac{a_k}{a_k+b_k}\right)^{\lambda_k} \le \sum\limits_{k=1}^n \lambda_k\left(\frac{a_k}{a_k+b_k}\right),\end{eqnarray*}
Similarly, we get
\begin{eqnarray*}\label{e10}\prod\limits_{k=1}^n\left(\frac{b_k}{a_k+b_k}\right)^{\lambda_k} \le \sum\limits_{k=1}^n \lambda_k\left(\frac{b_k}{a_k+b_k}\right).\end{eqnarray*}
Summing up these two inequalities, we have
\begin{eqnarray*}\prod\limits_{k=1}^n \frac{1}{\left(a_k+b_k \right)^{\lambda_k}} \left[ \prod\limits_{k=1}^n a_k^{\lambda_k}+\prod\limits_{k=1}^n b_k^{\lambda_k}\right] \le \sum\limits_{k=1}^n \lambda_k =1,\end{eqnarray*}
which leads to the desired result (\ref{e8}).

\begin{remark} A particular case $b_1=b_2=\cdots=b_n=1,\lambda_1=\lambda_2=\cdots=\lambda_n=\frac{1}{n}$ in (\ref{e8}) yields
\begin{eqnarray*} (1+a_1)(1+a_2)\cdots(1+a_n) \ge \left[1+\left(a_1a_2 \cdots a_n \right)^{\frac{1}{n}} \right]^n,\end{eqnarray*}
which is a famous inequality, called Chrystal inequality(see\cite[pp. 61]{GJG}), so we can view lemma \ref{lemma2.1} as a generalization of Chrystal inequality.
\end{remark}

\noindent
{\bf Proof 3.}
Use induction on $n$. When $n=1$, the result is obvious. Assume that (\ref{e5}) is true for $n=m$, that is
\begin{eqnarray*} \frac{a_1^r}{b_1^s}+\frac{a_2^r}{b_2^s}+\cdots+\frac{a_m^r}{b_m^s} \ge \frac{\left(a_1+a_2+\cdots+a_m\right)^r}{m^{r-s-1}\left(b_1+b_2+\cdots+b_m\right)^s}.\end{eqnarray*}
When $n=m+1$,
we need to prove the following inequality:
\begin{eqnarray*}\begin{aligned}\sum\limits_{k=1}^{m+1} \frac{a_k^r}{b_k^s} &= \sum\limits_{k=1}^m \frac{a_k^r}{b_k^s} +\frac{a_{m+1}^r}{b_{m+1}^s}\\
&\ge \frac{\left(a_1+a_2+\cdots+a_m\right)^r}{m^{r-s-1}\left(b_1+b_2+\cdots+b_m\right)^s}+\frac{a_{m+1}^r}{b_{m+1}^s}~~~\hbox{(by induction assumption)}\\
&= \frac{\left[ \Bigl(R_m(a)+\frac{a_{m+1}^r}{b_{m+1}^s}\Bigr)^{\frac{1}{r}}\bigl(S_m(b)+b_{m+1}\bigr)^{\frac{s}{r}} \left(m+1\right)^{\frac{r-s-1}{r}} \right]^r}{(m+1)^{r-s-1}(S_m(b)+b_{m+1})^s}\\
&\ge \frac{\left[\bigl(R_m(a)\bigr)^{\frac{1}{r}}\bigl(S_m(b)\bigr)^{\frac{s}{r}}m^{\frac{r-s-1}{r}}+\bigl(\frac{a_{m+1}^r}{b_{m+1}^s}\bigr)^{\frac{1}{r}}b_{m+1}^{\frac{s}{r}}1^{\frac{r-s-1}{r}}\right]^r}{(m+1)^{r-s-1}(b_1+\cdots+b_m+b_{m+1})^s}
\\ & ~~~~~~~~~~~~~~~~~~~~~~~~~~~~~~~~~~~~~~~~~~~\hbox{(by a special case $n=3$ in (\ref{e8}))}\\
&= \frac{(a_1+\cdots+a_m+a_{m+1})^r}{(m+1)^{r-s-1}(b_1+\cdots+b_m+b_{m+1})^s},\end{aligned}\end{eqnarray*}
where $R_m(a)=\frac{(a_1+\cdots+a_m)^r}{m^{r-s-1}(b_1+\cdots+b_m)^s}$, $S_m(b)=b_1+b_2+\cdots+b_m$.
Thus, inequality (\ref{e5}) holds for $n=m+1$, so the proof of the induction step is complete.

In the next theorem, we will prove equivalence relation between the weighted power mean inequality and Radon inequality, which is partly motivated by a slight observation of inequality (\ref{e7}).

\begin{theorem} The Radon inequality \eqref{e2} is equivalent to the weighted power mean inequality \eqref{e3}.\end{theorem}

\noindent
{\bf Proof .} $\Longrightarrow$ By the Radon inequality \eqref{e2} and $y_1,y_2,\ldots,y_n \in [0,+\infty),$ we have
\begin{eqnarray*}\begin{aligned}p_1y_1^{\frac{r}{s}}+p_2y_2^{\frac{r}{s}}+\cdots+p_ny_n^{\frac{r}{s}} &= \frac{(p_1y_1)^{\frac{r}{s}}}{p_1^{\frac{r}{s}-1}}+\frac{(p_2y_2)^{\frac{r}{s}}}{p_2^{\frac{r}{s}-1}}+\cdots+\frac{(p_ny_n)^{\frac{r}{s}}}{p_n^{\frac{r}{s}-1}} \\
&\ge \frac{ \left(p_1y_1+p_2y_2+  \cdots +p_ny_n \right)^{\frac{r}{s}}}{\left( p_1+p_2+ \cdots +p_n \right)^{\frac{r}{s}-1}}.\end{aligned}\end{eqnarray*}
which means that
\begin{eqnarray}\label{e4} \frac{p_1y_1^{\frac{r}{s}}+p_2y_2^{\frac{r}{s}}+\cdots+p_ny_n^{\frac{r}{s}}}{p_1+p_2+\cdots+p_n} \ge \left( \frac{p_1y_1+p_2y_2+\cdots+p_ny_n}{p_1+p_2+\cdots+p_n} \right)^{\frac{r}{s}}.\end{eqnarray}
Let $y_k=x_k^s$ for all  $x_k \ge 0$ $(k=1,2,\ldots,n)$ in (\ref{e4}). Thus, we can obtain the following weighted power mean inequality \eqref{e3}
\begin{eqnarray*} \left( \frac{p_1x_1^r+p_2x_2^r+\cdots+p_nx_n^r}{p_1+p_2+\cdots+p_n}\right)^{\frac{1}{r}} \ge \left( \frac{p_1x_1^s+p_2x_2^s+\cdots+p_nx_n^s}{p_1+p_2+\cdots+p_n}\right)^{\frac{1}{s}}.\end{eqnarray*}
$\Longleftarrow$ Let $p_k=b_k,x_k=\frac{a_k}{b_k}$ and $r=m+1(m\ge0), s=1$ in (\ref{e3}). Then, we have
\begin{eqnarray*} \left[ \frac{1}{b_1+b_2+\cdots+b_n} \left( \frac{a_1^{m+1}}{b_1^m}+ \frac{a_2^{m+1}}{b_2^m}+\cdots+ \frac{a_n^{m+1}}{b_n^m}\right) \right]^{\frac{1}{m+1}} \ge \frac{a_1+a_2+\cdots+a_n}{b_1+b_2+\cdots+b_n},\end{eqnarray*}
which implies that the Radon inequality (\ref{e2}) is achieved.

\begin{theorem}The following inequalities are equivalent:\\
(i) Bernoulli inequality,\\
(ii) the weighted AM-GM inequality,\\
(iii) H\"{o}lder inequality,\\
(iv) the weighted power mean inequality,\\
(v) Minkovski inequality,\\
(vi) Radon inequality. \end{theorem}

\noindent
{\bf Proof .}The equivalence between (iv) and (vi) is given in Theorem 2.2, the equivalence among (i), (iii) and (vi), one can find in \cite{DMB} as well as (ii), (iii) and (iv) in \cite{YTL}, the equivalence between (iii) and (v) is shown in \cite{LECH}.

\begin{corollary} If $a_1, a_2,\ldots,a_n,b_1,b_2,\ldots,b_n$ are positive real numbers,  then for $m\le-1$, the following inequality holds
\begin{eqnarray}\label{e11} \frac{a_1^{m+1}}{b_1^m}+ \frac{a_2^{m+1}}{b_2^m}+\cdots+ \frac{a_n^{m+1}}{b_n^m} \ge \frac{\left(a_1+a_2+\cdots+a_n \right)^{m+1}}{ \left(b_1+b_2+\cdots+b_n \right)^m}.\end{eqnarray} \end{corollary}

\noindent
{\bf Proof .} Since $m \le -1$, thus by the inequality (\ref{e2}), we have
 \begin{eqnarray*}\begin{aligned} \frac{a_1^{m+1}}{b_1^m}+ \frac{a_2^{m+1}}{b_2^m}+\cdots+ \frac{a_n^{m+1}}{b_n^m} &=\frac{b_1^{-m}}{a_1^{-m-1}}+ \frac{b_2^{-m}}{a_2^{-m-1}}+\cdots+ \frac{b_n^{-m}}{a_n^{-m-1}}\\
&\ge \frac{ \left(b_1+b_2+\cdots+b_n \right)^{-m}}{ \left(a_1+a_2+\cdots+a_n \right)^{-m-1}}.\end{aligned}\end{eqnarray*}
The inequality (\ref{e11}) holds.

\begin{corollary} If $a_1, a_2,\ldots,a_n,b_1,b_2,\ldots,b_n$ are positive real numbers, then for nonpositive real numbers $r, s$  such that $r\ge s+1$ , we have
\begin{eqnarray}\label{e12}\frac{a_1^r}{b_1^s}+\frac{a_2^r}{b_2^s}+\cdots+\frac{a_n^r}{b_n^s} \ge \frac{\left(a_1+a_2+\cdots+a_n\right)^r}{n^{r-s-1}\left(b_1+b_2+\cdots+b_n\right)^s}.\end{eqnarray} \end{corollary}

\noindent
{\bf Proof .} For $r\le 0, s\le 0$,  the inequalities $-s\ge -r+1,-r\ge 0,-s\ge 0 $ hold. By the inequality (\ref{e5}), we obtain
\begin{equation*}\begin{aligned} \frac{a_1^r}{b_1^s}+ \frac{a_2^r}{b_2^s}+\cdots+ \frac{a_n^r}{b_n^s} &= \frac{b_1^{-s}}{a_1^{-r}}+ \frac{b_2^{-s}}{a_2^{-r}}+\cdots+ \frac{b_n^{-s}}{a_n^{-r}} \\
&\ge \frac{ \left(b_1+b_2+\cdots+b_n \right)^{-s}}{n^{-s-(-r)-1} \left(a_1+a_2+\cdots+a_n \right)^{-r}}\\
&=\frac{\left(a_1+a_2+\cdots+a_n\right)^r}{n^{r-s-1}\left(b_1+b_2+\cdots+b_n\right)^s}.\end{aligned}\end{equation*}
So, the inequality (\ref{e12}) holds.

\begin{corollary} If $a_1, a_2,\ldots,a_n,c_1,c_2,\ldots,c_n$ are positive real numbers,  and $m$ is real numbers such that $m>0$ or $m\le -1$, then
\begin{equation}\label{e13} \frac{a_1}{c_1}+\frac{a_2}{c_2}+\cdots+\frac{a_n}{c_n} \ge \frac{(a_1+a_2+\cdots+a_n)^{m+1}}{\left(a_1c_1^{\frac{1}{m}}+a_2c_2^{\frac{1}{m}}+\cdots+a_nc_n^{\frac{1}{m}} \right)^m}.\end{equation}\end{corollary}

\noindent
{\bf Proof .} Consider $b_k=a_kc_k^{\frac{1}{m}}$ for all $1\le k \le n$  in the inequality (\ref{e2}) and (\ref{e11}).  Thus, we obtain the inequality (\ref{e13}).

\begin{corollary} If $a,b\in \mathbb{R},a<b,m\ge 0$ or $m\le-1$,$f,g:[a,b]\rightarrow(0,+\infty) $ are integrable functions on $[a,b]$ for any $x\in[a,b]$, then
\begin{equation}\label{e14} \int_a^b \frac{\left(f(x)\right)^{m+1}}{\left(g(x)\right)^m}dx \ge \frac{\left(\int_a^bf(x)dx\right)^{m+1}}{\left(\int_a^bg(x)dx\right)^m}.\end{equation}
\end{corollary}

\noindent
{\bf Proof .} Letting $n \in \mathbb{N}_+, x_k=a+k \frac{b-a}{n},k \in \{0,1,\ldots,n\}$ and $\xi_k \in [x_{k-1}, x_k]$. By inequality (\ref{e2}) and (\ref{e11}), we get
 \begin{equation*} \sum\limits_{k=1}^n \frac{\left(f(\xi_k)\right)^{m+1}}{\left(g(\xi_k)\right)^m} \ge  \frac{\left( \sum\limits_{k=1}^n f(\xi_k) \right)^{m+1}}{\left(\sum\limits_{k=1}^n g(\xi_k) \right)^m}.\end{equation*}
It results that
\begin{equation*} \sigma \left( \frac{(f(x))^{m+1}}{(g(x))^m}, \Delta_n, \xi_k \right) \ge \frac{\bigl[\sigma \left( f(x), \Delta_n, \xi_k \right)\bigr]^{m+1}}{\bigl[\sigma \left(g(x),\Delta_n,\xi_k \right)\bigr]^m},\end{equation*}
where $\sigma \left( f(x), \Delta_n, \xi_k \right)$ is the corresponding Riemann sum of function $f(x)$, of $ \Delta_n=(x_0,x_1,\ldots,x_n)$ division and the intermediate $\xi_k$ points. By passing to limit  in ineuqality above, when $n$ tends to infinity, the inequality(\ref{e14}) follows.

\begin{corollary} If $a,b\in {\bf R},a<b,rs\ge0,r \ge s+1$, $f,g:[a,b]\rightarrow(0,+\infty) $ are integrable functions on $[a,b]$ for any $x\in[a,b]$, then
\begin{equation}\label{e15} \int_a^b \frac{\left(f(x)\right)^r}{\left(g(x)\right)^s}dx \ge \frac{\left(\int_a^bf(x)dx\right)^r}{(b-a)^{r-s-1}\left(\int_a^bg(x)dx\right)^s}.\end{equation}
\end{corollary}

\begin{proposition}Show that if $a,b,c$ are the lengths of the sides of a triangle and  $2S=a+b+c$, then
\begin{equation}\label{m1} \frac{a^n}{b+c}+\frac{b^n}{c+a}+\frac{c^n}{a+b} \ge \left( \frac{2}{3} \right
)^{n-2} S^{n-1},\quad n \ge 1.\end{equation}\end{proposition}

\noindent
{\bf Proof .} When $n=1$, the result \eqref{m1} is Nesbitt inequality(see \cite[p. 16, example 1.4.8]{RJR} or \cite[p. 2, exercise 1.3]{ZC}).
For $n \ge 2$, by  (\ref{e5}), we have
\begin{eqnarray*} \frac{a^n}{b+c}+\frac{b^n}{c+a}+\frac{c^n}{a+b} &\ge& \frac{(a+b+c)^n}{3^{n-1-1}(b+c+c+a+a+b)}\\ &=& \left( \frac{2}{3} \right)^{n-2} S^{n-1}.\end{eqnarray*}

\begin{proposition} Let $a_1,a_2,\ldots,a_n$ be positive real numbers such that $a_1+a_2+\cdots+a_n=s$ and $p \geqslant q+1 \geqslant 1$. Prove that
\begin{equation*} \sum\limits_{k=1}^n \frac{a_k^p}{(s-a_k)^q} \ge \frac{s^{p-q}}{(n-1)^q n^{p-q-1}}.\end{equation*}\end{proposition}

\noindent
{\bf Proof .} By applying the inequality (\ref{e5}), the inequality above is easily obtained.

\begin{proposition}Let $x,y$, and $z$ be positive real numbers such that $xyz=1$. Then
\begin{equation*} \frac{x^3}{(1+y)(1+z)}+ \frac{y^3}{(1+z)(1+x)}+\frac{z^3}{(1+x)(1+y)} \ge \frac{3}{4}.\end{equation*}\end{proposition}

\noindent
{\bf Proof .} By the generalized Radon inequality (\ref{e5}), we obtain
\begin{equation*} \begin{aligned} & \frac{x^3}{(1+y)(1+z)}+ \frac{y^3}{(1+z)(1+x)}+\frac{z^3}{(1+x)(1+y)}\\ &\ge \frac{(x+y+z)^3}{3 \left( (1+y)(1+z)+(1+z)(1+x)+(1+x)(1+y) \right)}\\ &= \frac{(x+y+z)^3}{ 9+6(x+y+z)+3(xy+yz+zx)} \\ & ~~~~~~~~~~~~\hbox{(by a general inequality $3(xy+yz+zx) \le (x+y+z)^2$ )}\\ &\ge \frac{(x+y+z)^3}{ 9+6(x+y+z)+(x+y+z)^2}.\end{aligned} \end{equation*}
Since $x+y+z \ge 3 \sqrt[3]{xyz}=3$, it is easy to prove that $\frac{(x+y+z)^3}{ 9+6(x+y+z)+(x+y+z)^2} \ge \frac{3}{4}$. Another proof can be found in \cite[pp. 139-140]{RJR}.

\vspace{1cm}\noindent{\textbf{References}}
\bibliographystyle{elsarticle-num}
\bibliography{<your-bib-database>}

\begin{thebibliography}{00}
    \normalsize
\baselineskip=17pt

\bibitem{KYF} K.Y. Fan, \textit{Generalization of Bergsrr\"{o}m inequality}, Amer. Math.
Monthly, \textbf{66} (1959), 153--154.

\bibitem{HB} H. Bergstr\"{o}m, \textit{A triangle inequality for matrices}, Den Elfte Skandinaviske
Matematikerkongress, Trondheim, 1949, Johan Grundt Tanums Forlag, Oslo, 1952, 264--267.

\bibitem{RB} R. Bellman, \textit{Notes on matrix theory-IV (An inequality due to Bergstr\"{o}m)},
Amer. Math. Monthly, \textbf{62} (1955), 172--173.

\bibitem{SBJ} S. Abramovich, B. Mond, J. Pe\v{c}ari\'{c}, \textit{Sharpening Jensen's inequality
and a majorization theorem}, J. Inequal. Pure Appl. Math., \textbf{214} (1997), 721--728.

\bibitem{JFY} J.E. Pe\v{c}ari\'{c}, F. Proschan, Y.I. Tong, \textit{Convex functions, partial
ordering, and statistical operations}, Academic Press, San Diego, 1992.

\bibitem{JR} J. Radon, \textit{\"{U}ber die absolut additiven Mengenfunktionen}, Wiener
Sitzungsber \textbf{122} (1913), 1295--1438.

\bibitem{GJG} G. Hardy, J.E. Littlewood, G. Polya, \textit{Inequalities}, Cambridge University
Press, Cambridge, 1934.

\bibitem{CBM} C.B. Morrey, \textit{A class of representations of manifolds}, Amer. J. Math.,
\textbf{55} (1933), 683--707.

\bibitem{RJR} R.B. Manfrino, J.A.G. Ortega, R.V. Delgado, \textit{Inequalities: A Mathematical
Olympiad Approach}, Birkh\"{a}user, Basel-Boston-Berlin, 2009.

\bibitem{CM} C. Mortici, \textit{A new refinement of the Radon inequality}, Math. Commun.,
\textbf{16} (2011), 319-324.

\bibitem{DMB}D.M. B\v{a}tinetu-Giurgiu and O.T. Pop, \textit{A generalization of Radon's
Inequality}, Creative Math. and Inf., \textbf{19} (2) (2010), 116-121.

\bibitem{ZC} Z. Cvetkovski, \textit{Inequalities. Theorems, techniques and selected problems},
Springer, Heidelberg, 2012.

\bibitem{DJA} D.S. Mitrinovi\'{c}, J.E. Pe\v{c}ari\'{c} and A.M. Fink, \textit{Classical and
new inequalities in analysis. Mathematics and its applications (East European Series)}, 61.
Kluwer Amcademic Publishers Group, Dordrecht, 1993.

\bibitem{YKC} K.-C. Yang, \textit{A note and generalization of a fractional inequality},
Journal of Yueyang Normal University (Natural Sciences Edition), \textbf{15} (4)
(2002), 9-11. (In Chinese) DOI: \href{http://dx.chinadoi.cn/10.3969\%2fj.issn.1672-5298.2002.04.003}
{10.3969/j.issn.1672-5298.2002.04.003}.

\bibitem{YTL}Y. Li and X.-M. Gu, \textit{The weighted AM-GM inequality is equivalent to the H\"{o}lder inequality},
submitted to Anslysis, in revision, November 24, 2015, 6 pages. Also available online
at \url{https://arxiv.org/abs/1504.02718}.

\bibitem{LECH}
L. Maligranda, \textit{Equivalence of the H\"{o}lder-Rogers and Minkowski Inequalities}, Math. Inequal. Appl., \textbf{4} (2011), 203-207.
\end{thebibliography}
    
\end{document}